\documentstyle{amsart}
\date{  \today}
\def\appendix{\par
 \setcounter{section}{0}
 \setcounter{subsection}{0}
 \def\thesection{Appendix \Alph{section}}}

\newtheorem{th}{Theorem}

\newtheorem{dfn}{Definition}
\newtheorem{lem}{Lemma}

\newcommand{\rar}{\rightarrow}
\newcommand{\bfp}{{\Bbb{P}}}
\newcommand{\bfc}{{\Bbb{C}}}
\newcommand{\bfq}{{\Bbb{Q}}}
\newcommand{\bff}{{\Bbb{F}}}
\newcommand{\bfz}{{\Bbb{Z}}}

\newcommand{\clx}{{\cal{X}}}
\newcommand{\cle}{{\cal{E}}}
\newcommand{\co}{{\cal{O}}}
\newcommand{\calp}{{\cal{P}}}
\newtheorem{prp}{Proposition}

\newcommand{\spec}{{\mbox{Spec }}}
\newcommand{\sym}{{\mbox{Sym}}}
\newtheorem{cor}{Corollary}

\newcommand{\oks}{{{\cal O}_{K,S}}}

\begin{document}
\begin{center}{\vspace*{.3cm}}{ \bf \LARGE Uniformity of stably integral
points on elliptic curves}\\[.5cm] by
\\ {\sc Dan  Abramovich}\footnote{Partially supported by NSF grants
 DMS-9207285 and DMS-9503276}\\  \today
\end{center}

\addtocounter{section}{-1}
\section{Introduction}

Let $X$ be a variety of logarithmic general type, defined over a number field
$K$. Let $S$ be a finite set of places in $K$ and let $\oks$ be the ring of
$S$-integers. Suppose that $\clx$ is a model of $X$ over $\spec\oks$.
 As a natural generalizasion of theorems of Siegel and Faltings,  It was
conjectured by S. Lang and P. Vojta (\cite{vojta}, conjecture 4.4) that the
set of $S$-integral points $\clx(\oks)$ is not
Zariski dense in $\clx$. In case $X$ is projective, one may chose an arbitrary
projective model $\clx$ and then $\clx(\oks)$ is identified with $X(K)$. In
such a case, one often refers to this conjecture of Lang and Vojta as just
Lang's conjecture.

L. Caporaso, J. Harris
and B. Mazur \cite{chm} apply Lang's conjecture in the following way: Let
$X\rar B$
be a smooth family of curves of genus $g>1$. Let $X^n_B\rar B$ be the $n$-th
fibered
power of $X$ over $B$. In  \cite{chm} it is shown that for high enough $n$,
the variety $X^n_B$ dominates a variety of general type. Assuming Lang's
conjecture, they deduce the following remarkable result: the number of rational
points on a curve of genus $g$ over a fixed number field is uniformly bounded.

In this note we study an analogous implication for elliptic curves. Let $E/K$
be an elliptic curve over a number field, and let $P\in E(K)$. We say that $P$
is stably $S$-integral, denoted $P\in E(K,S)$, if $P$ is $S$-integral after
semistable reduction (see \S\ref{stably}). Our main theorem states (see
\S\ref{main}):
\begin{th}(Main theorem in terms of points) Assume that the Lang - Vojta
conjecture
holds. Then for any number field $K$ and a finite set of places $S$, there is
an
integer $N$ such that for any elliptic curve $E/K$ we have $\#E(K,S)<N$.
\end{th}

Since the moduli space of elliptic curves is only one-dimensional, the
computations and the proofs are a bit  simpler than the higher genus cases.
One can view the results in this paper as a simple application of
the methods of \cite{chm}.

\subsection{Overview} In section \ref{correlation} we prove a basic lemma
analogous to lemma 1.1 \cite{chm} on uniformity of correlated points. In
section \ref{level3} we study a particular pencil of elliptic
curves which is the main building block for proving theorem 1. In section
\ref{twist} we look at quadratic twists of an elliptic curve, motivating the
study in section
\ref{stably} of stably integral points. Section \ref{main} gives a proof of the
main theorem.

In section \ref{ubc} we will refine our methods and show that the Lang - Vojta
conjecture implies the uniform boundedness conjecture for torsion on elliptic
curves, thus giving a conditional (and therefore obsolete) proof of the
following theorem of Merel:
\begin{th}(Merel, \cite{merel}) For and integer $d$ there is an integer $N(d)$
such that given a number field $K$ with $[K:\bfq]=d$, and given an elliptic
curve $E/K$, then $\#E(K)_{tors}<N(d)$.
\end{th}
It should be noted that the methods introduced in section \ref{ubc} were
essential for the developments in \cite{abr},\cite{abr1}.
It would be interesting if one could apply the method to the study
points on abelian varieties in general.

\subsection{Acknowledgements} I am indebted to  the ideas in the
work \cite{chm} of Lucia Caporaso, Joe Harris and Barry Mazur and to
conversations with
them. Thanks to Joe Silverman and Sinnou David for discussions of points on
elliptic curves and much encouragement. Much of this text was written during a
pleasant
visit with the group ``Problemes Diophantiens'' in Paris. It is also a pleasure
to thank Henri Darmon, Gerhard Frey Brendan Hassett and Felipe Voloch for
their suggestions. Special thanks are due to Frans Oort who read an
earlier version of this paper and sent me many helpful comments\footnote{I
take the opportunity to wish Professor Oort a happy 60th birthday.}.

\subsection{The Lang - Vojta conjecture}\label{lv}
A common practice in arithmetic geometry is that of generalizing rational
points on projective varieties to integral points on quasi-projective
varieties. We can summarize this in the following table, which will be
explained
below:\\

\begin{tabular}{l|l}\hline\hline
Number field $K$ & Ring of $S$-integers $O_{K,S}$ \\[2mm] \hline
Projective variety $X$ over $K$ & \parbox{3in}{\vspace*{1mm} Quasi projective
variety $X$
and a model
$\clx$ over $O_{K,S}$ } \\ \hline \vspace*{1mm}
Rational point $P\in X(K)$ & Integral point $P\in \clx(\oks)$ \\[2mm] \hline
$X$ of general type & $X$ of log-general type \\
e.g.: $C$ a curve of genus $>1$ & \parbox{3in}{ e.g.: $E$ an elliptic curve
with the origin
removed } \\ \hline
Faltings' theorem: $C(K)$ finite & Siegel's theorem: $E(O_{K,S})$ finite
\\
\hline \vspace*{1mm}
 \parbox{3in}{\vspace*{1mm} Lang's conjecture: If $X$ is of general type then
$X(K)$ not
Zariski
dense\\ } &
 \parbox{3in}{\vspace*{1mm} Lang-Vojta conjecture: If  $X$ is of logarithmic
general type
then
$\clx(\oks)$ not
Zariski dense } \\ \hline
 \parbox{3in}{\cite{chm}: Lang's conjecture implies uniformity of $\#C(K)$} &
\parbox{3in}{
\vspace*{1mm} ?? }\\ \hline\hline
\end{tabular}
\vspace*{2mm}

We remind the reader of the definition of a variety of  log general type:
\begin{dfn} Let $X$ be a quasi-projective variety over $\bfc$. Let $f:Y \rar X$
be a resolution of singularities, that is, a proper, birational morphism where
$Y$ is a smooth variety. Let $Y\subset Y_1$ be a projective compactification,
such that $Y_1$ is smooth and such that
$D = Y_1\setminus Y$ is a divisor of normal crossings. Then $X$ is said to be
of logarithmic general type if for some positive integer $m$, the rational map
defined by the complete linear system $|m(K_{Y_1}+D)|$ is birational to the
image.
\end{dfn}

 Let $X$ be a quasi-projective variety
of logarithmic general type, defined over a field $K$ which is finitely
generated over $\bfq$ (e.g., a number field). Let $R$ be a ring, finitely
generated over $\bfz$, whose
fraction field is $K$ (e.g., the ring of $S$-integers in a number field).
Choose a model $\clx$ of $X$ over $R$. The following is
a well-known
conjecture of Lang and Vojta (\cite{vojta}, conjecture 4.4).

{\bf Conjecture.} The set of integral points $\clx(R)$ is not Zariski dense in
$\clx$.

In case $X$ in the conjecture above is projective, then logarithmic general
type means just general type; and integral points are just rational points.

\subsection{What should the last entry in the table read?}
We would like to fill in the question mark in the last entry in the table. One
is tempted to ask:

{\em Does the Lang - Vojta  conjecture imply the uniformity of
$\#E(O_{K,S})$?}\\
but one sees immediately that this cannot be true without some restrictions.
Most importantly, one has to restrict the choice of the model $E$, as can be
seen in the following example:

Let $E$ be an elliptic curve over a number field $K$ such that $E(K) $ is
infinite. Fix $P_1,\ldots ,P_n\in E(K)$. Choose an equation $$y^2 = x^3 +
Ax+B$$ for $E$, where $A,B\in O_{K,S}$ for some finite $S$. Choose $c\in \oks$
such that for each $P_i$ one has $c^2x(P_i),c^3 y(P_i)\in \oks$. By changing
coordinates $x_1 = c^2 x, y_1 = c^3y$, one obtains a new model $E_1$ given by
the equation $y_1^2 = x_1^3 + c^4 A x_1 + c^6B$, on which all the points $P_i$
are integral.

The problem with this new model arises because when one changes coordinates,
one blows up the closed point corresponding to the origin at primes dividing
$c$, so the resulting model has ``extraneous'' components over these primes. We
are led to modify the statement:

\begin{th}(Main theorem in terms of models, see \S\ref{stably}) Assume that the
Lang-Vojta
conjecture holds. Then for any number field $K$ and a finite set of places $S$
there is an integer  $N$ such that  for any stably
minimal elliptic curve $\cle$ over $\oks$  we have $\#\cle(\oks)<N$.
\end{th}

It turns out that stably minimal models are very minimal indeed. In
particular we will see that N\'eron models, the  canonical models of
elliptic curves over rings of integers, are not necessarily sufficiently
minimal for the purpose of
our methods. On the other hand
we will see that semistable models are stably minimal. A precise
definition of a stablly minimal model, and how to obtain a canonical one
from the N\'eron model, will be given in \S\ref{stably}.

 In the case of  semistable elliptic curves,
it is worthwhile to state an immediate corollary of the theorem:

\begin{cor} The Lang-Vojta conjecture implies that the number of
integral points on  semistable elliptic curves over $\bfq$ is bounded.\end{cor}

{\bf Remark.} Another conjecture of Lang (see \cite{lang1}) predicts that the
number of all
$S$-integral points on so called {\em quasi-minimal} elliptic curves  should be
bounded in
terms of the rank and the number of elements in $S$. In view of the corollary,
one is tempted to ask whether the rank of an elliptic curve can be bounded in
terms of the places of additive reduction of the elliptic curve.

\section{Boundedness of correlated points}\label{correlation}

One of the main ideas in \cite{chm} is, that in order to bound the number of
points on curves it is enough to show that they are {\em correlated}, that is,
there is
an algebraic relation between all $n$-tuples of these points. This is the
content of the lemma below. First,
some notation.

Let $\pi:X\rar B$ be a family of smooth irreducile curves over a field $K$. We
denote by
$\pi_n:X_B^n\rar B$ the n-th
fibered power of $X$ over $B$. Given a point $b\in B$ we denote by $X_b$ the
fiber of $X$ over $b$; Similarly, given  $Q=(P_1,\ldots,P_n)\in X^n_B$ we
denote by $X_Q\subset X_B^{n+1}$ the fiber of $X_B^{n+1}$ over $Q$.  Note that
if $\pi_n(Q) = b$ then $X_Q \simeq X_b$. Denote by $p_n:X_n\rar X_{n-1}$ the
projection onto the first $n-1$ factors.

Assume that we are given a
subset $\calp\subset X(K)$ (typical examples would be rational points, or
integral
points on some model of $X$). Again, we  denote by $\calp_B^n\subset
X_B^n$ the fibered power of $\calp$ over $B$ (namely the union of the
$n$-tuples of points in $\calp$ consisting of points in the same fiber), and by
$\calp_b$ the points of $\calp$ lying over $b$.

\begin{dfn} Assume
that for some $n$ there is  a proper closed subset $F_n\subset
X^n_B$ such that $\calp_B^n\subset F_n$. In such a case we say that the
subset $\calp$ is {\em $n$-correlated}.
\end{dfn}

 For instance, a subset $\calp$ is
1-correlated if and only if it is not Zariski dense; in which case it is easy
to see that, over some open set in $B$, the number of points of $\calp$ in each
fiber is bounded. This is generalized by the following lemma:

\begin{lem}(compare \cite{chm}, lemma 1.1) Let $X\rar B$ be a family of smooth
irreducible curves, and let
$\calp\subset
X(K)$ be an $n$-correlated subset.
 Then there is a dense open
set $U\subset B$ and an integer $N$ such that for every $b\in U$, we have
 $\# \calp_b\leq N$.
\end{lem}

{\bf Proof:} Let $F_n = \overline{\calp_B^n}$ be the Zariski closure, and
$U_n= X_B^n\setminus F_n$ the complement. We now define be descending
induction: $U_{i-1}  =
p_i (U_i)$ and $F_{i-1}=X^{i-1}_B\setminus U_{i-1} $ the complement. Notice
that
over $U_{i-1}$, the map $p_i$ restricts to a finite map on $F_i$: by definition
if $x\in U_{i-1}$ then $p_i^{-1}(x)\not\subset F_i$, and $p_i^{-1}(x) $ is
an irreducible curve. Therefore
the number of points in the fibers of this map is bounded: if $x\in U_{i-1}$
then we can write $\#(p_i^{-1}(x)\cap F )\leq d_i$.

Let $U=U_0\subset B$. We claim that over $U$, the number of  points of $\calp$
in each fiber is  bounded. Consider a point $b\in U$.

Case 1:  $\calp_b \subset F_1$. In this case, the number of points on
$\calp_b$ is bounded by $d_1$.

Case 2: there is some $P\in  \calp_b, P\not\in F_1$, but $X_P\cap\calp_b^2
\subset F_2$.
In this case the number of points is bounded by $d_2$.

Case $i$: $Q=(P_1,\ldots,P_{i-1})\in \calp_b^{i-1}\setminus F_{i-1}$ but
$X_Q\cap \calp_b^i\subset
F_i$. Here the number of points is bounded
by $d_i$.

Notice that in the case $i=n$ we  have  by definition $X_Q\cap \calp_B^n\in
F_n$, and the process stops. Therefore $N=\displaystyle \max_i d_i$ is a bound
for the number of $\calp$ points in each fiber over $U$.

{\bf Example} (\cite{chm}): Let $X\rar B$ be a family of smooth, irreducible
curves of
genus $>1$ over a number field $K$. Assume that Lang's conjecture holds true.
Then in \cite{chm} it is shown that $X(K)$ is $n$-correlated, and the lemma
above, with noetherian induction, is used to obtain the existence of a uniform
bound on the number of rational points on such curves.

{\bf Example:} Assume that $X_K\rar B_K$ is a semistable family of curves of
genus 1,
together
with a section $s:B_K\rar X_K$, and assume
that over an open set $B_0\subset B_K$ the restricted family $X_0\rar B_0$  is
smooth. Assume that the Lang - Vojta conjecture holds true.  Given a semistable
model $X$ of $X_K\setminus s(B_K)$ over $\oks$, we
will later show that $X(\oks)\cap X_0$ is $n$-correlated for some integer $n$.
We will deduce the existence of a uniform bound on the number of integral
points on curves in this family.

{\bf Example:} Let $S$ be  a finite set of places in $K$. Assuming that the
Lang - Vojta
conjecture holds true, we will show that the set of
stably $S$-integral points on any family of elliptic curves over a number field
$K$ is $n$-correlated for some $n$. We will deduce the existence of a uniform
bound on the number of stably $S$-integral points on an elliptic curve.

\section{Moduli of elliptic curves with level 3 structure}\label{level3} We
introduce here a building block in the proof of the main theorem.
\subsection{The geometry} Let $E_1$ be the universal family of elliptic
curves over $\bfc$ with full symplectic level 3 structure. The surface $E_1$
can be identified with the total space of the elliptic
pencil written in bi-homogeneous coordinates as: $$(*)\quad
\lambda(X^3+Y^3+Z^3)
- 3\mu XYZ=0,$$ mapping to the moduli space $\bfp^1$ via
$[\lambda:\mu]$. This equation gives a smooth model, which by abuse of notation
we will also call $E_1$, of this space over   $\spec\bfz[1/3]$. We
may choose the section $\Theta$ over the point $[X:Y:Z]=[1:-1:0]$ as the origin
of
the elliptic surface. Over $\spec \bfz[1/3]$, the fibers of the elliptic pencil
possess level 3 structure of type $\mu_3 \times \bfz/3\bfz$, in a way which is
described precisely by Rubin and Silverberg in \cite{rs}; however in this
section we will work over $\bfc$.

The pencil $E_1\rar \bfp^1$ is semistable, possessing four singular fibers,
having 3 nodes each, over $\Sigma_0=\{0,1,\zeta_3,\zeta_3^2\}\subset \bfp^1$
where $\zeta_3$ is a primitive third root of 1.

Let $L$ be the pullback of a line from the plane, and
let $S_1,\ldots,S_9$ be the exceptional curves over the nine base points of the
pencil fixing $S_1=\Theta$ to be the origin of the elliptic surface. Let
$F$ be a fiber of the elliptic surface. We have the linear equivalence  $-F\sim
-3L +S_1+\cdots+S_9$.

As a pencil of cubics with smooth total space, one easily calculates the
relative dualizing sheaf, as follows:
 We know that $\omega_{\bfp^2} =
\co_{\bfp^2}(-3)$. The canonical sheaf of the blown up surface is therefore
$\co( -3L +S_1+\cdots+S_9)$. Therefore we have
 $\omega_{E_1} \simeq \co(-F)$, and $\omega_{E_1/\bfp^1} \simeq \co(F)$.

Let $\pi_n:E_n\rar \bfp^1$ be
the $n$-th fibered power of $E_1$ over $\bfp^1$. Denote by $\pi_{n,i}:E_n\rar
E_1$ the projection onto the $i$-th factor. We have that
$\omega_{E_n/\bfp^1} \simeq \co(nF)$, and therefore $\omega_{E_n}
\simeq \co((n-2)F)$.  We denote by
$\Theta_n=\sum_{i=1}^n\pi_{n,i}^*\Theta$, the theta divisor. We denote by
$\Sigma_n=\pi_n^{-1}\Sigma_0\subset E_n$
the locus of singular
fibers, the inverse image of $\Sigma_0\subset \bfp^1$.

It should be noted that $E_n$ is singular, but not too singular:

 \begin{lem}
There is a desingularization  $f_n:\tilde{E_n}\rar{E_n}$ such that
$\omega_{\tilde E_n} \simeq f_n^*\omega_{ E_n}(D)$, for some effective divisor
$D$ such that $f_n(D)\subset \Sigma_n$, and such that $f_n^*\Theta_n$ is a
reduced divisor of normal crossings.
\end{lem}

{\bf Proof:} The existence of a desingularization with $\omega_{\tilde E_n}
\simeq f_n^*\omega_{ E_n}(D)$ follows from \cite{chm}, lemma 3.3, or lemma
3.6 of \cite{viehweg}. The desingularization is given by a succession of
blowups along smooth centers. Since the
singular locus of $E_n$ meets $\Theta_n$ transversally, the centers of the
blowups can be taken to be transversal to $\Theta_n$, and therefore its inverse
image is a divisor of normal crossings.

This lemma shows that $E_n\setminus\Theta_n$ has log canonical
singularities. This means that sections of powers of $\omega_{E_n}(\Theta_n)$
give regular sections of the logarithmic pluricanonical sheaves of
$\tilde{E_n}$; therefore, in order to prove that $E_n$ is of logarithmic
general type, there is no need to pass to a resolution of singularities - it
suffices to show that $\omega^k_{E_n}(k\Theta_n)$ has many sections.

\begin{lem} \begin{enumerate}
\item The line bundle $\omega_{E_1/\bfp^1}(\Theta_1)$ is the pullback of
an ample bundle along a birational morphism.
\item Fix  $n>2$. Then
 ${E_n}\setminus\Theta_n$ is of logarithmic general type. Moreover, the base
locus of the logarithmic pluricanonical linear series is contained in
$\Theta_n$.
\end{enumerate}
\end{lem}

{\bf Proof: } Let $Y$ be the blowup of $\bfp^2$ at all the base points of our
pencil except
$[1,-1,0]$. The surface $Y$ is the same as  $E_1$ blown down along $\Theta$.
The line
bundle $\omega_{E_1/\bfp^1}(\Theta_1)$ is the pullback of
$M=\co(3)\otimes\co(-(S_2+\cdots+S_9))$ from $Y$. On $Y$, $M$ is
represented by the strict transform of one of the cubics of the pencil, hence
it is a nef line bundle; it
has  self intersection number 1, therefore it is nef and big. In fact, it is
easy to see by a dimension count that the complete linear system of sections of
$M^{\otimes 3}$ gives a birational morphism of $Y$ to a surface in projective
space,
which blows down only the fibral components of $E_1$ which do not meet $S_1$.

Part (2)  follows by taking the products of sections pulled back along the
projections:
On $E_3$, let $p_{ij},p_k$ be the projections to the  $i$-th and $j$-th
factors,
respectively $k$-th  factor. We have the inclusion  $p_{12}^*\omega_{E_2}
\otimes p_3^*\omega_{E_1/\bfp^1}(\Theta_1) \subset\omega_{E_3}(\Theta_3)$.
The sections of a power of this subsheaf give a map which generically
separates between
points whose third factors are different. Part 1 of this lemma implies
that the base locus of these sections is
contained in the theta divisors. By repeating this
for the other two
projections, we find that sections of powers of $\omega_{E_3}(\Theta_3)$
generically separate points; in particular $E_3\setminus\Theta_3$ is of
logarithmic general type. Similarly, $E_n$ is of logarithmic general type for
$n\geq 3$.

\subsection{Boundedness of integral points on elliptic curves with level 3
structure}

Let $K$ be a number field containing $\bfq(\zeta_3)$, and let $E$ be an
elliptic curve with full symplectic level 3
structure over $K$. Let $R={\cal O}_K[1/3]$. The curve $E$ occurs as a fiber in
the
surface over a point in $\bfp^1(K)$, which automatically has semistable
reduction over $R$. Let $\cle$ be the semistable model. Given any three
$R$-integral
points $P_i$ on $\cle\setminus 0$, the point $(P_{1}, P_{2}, P_{3})$ gives
rise to
an $R$-integral point of the scheme $E_3\setminus\Theta_3$  (where by abuse of
notation, we use the model of $E_3$ over $R$ which is the fibered cube of the
given model $(*)$ of $E_1$).

Assume that the Lang - Vojta conjecture holds for the variety $E_3\setminus
\Theta_3$. Thus the Zariski closure $F$ of the set of integral
points $(E_3\setminus\Theta_3)(R)$ is
a proper subvariety of $E_3$; in other words, the set  $\calp=(E\setminus
\Theta)(R)$ is 3-correlated. By lemma 1,
there is  a dense open set $U
\in \bfp^1$ such that  the number of integral points of fibers over  $U(K)$ is
bounded. The complement of
$U$ is a finite number of points, therefore by Siegel's theorem there is a
bound on the number of integral points on these curves as well.

{\bf Open problem:} Show that the $S$-integral points on $E_3\setminus
\Theta_3$
are not Zariski dense.
\section{Quadratic twists of an elliptic curve}\label{twist}

As a ``complementary case'' to the last section we will discuss here a typical
case of isotrivial families of elliptic curves.  This is in direct analogy with
the exposition in \cite{chm}, \S\S 2.2. It will give us a good hint about the
type of
models of elliptic curves we need in order to obtain boundedness. A slightly
more general version of the example here will be used in the proof of the main
theorem.

Let $E: y^2 = x^3 + Ax + B$ be a fixed elliptic curve. We  denote $f(x) =
x^3 + Ax + B$. We assume that $A$ and
$B$ are relatively prime $S$-integers in a number field $L$, where $S$ is a
finite set of
places.
All the quadratic twists of the curve over $L$ can be written in the form:
 $$E_t: ty^2 = f(x)$$
where $t$ may be chosen $S$-integral.

We may form the family of  Kummer surfaces associated to $E_t$:
 $$K_t: t^2 z^2 = f(x_1) f(x_2). $$ We have a morphism of $K_t$ to
$K_1$ via $(x_1,x_2,z,t) \mapsto (x_1,x_2,tz)$. It can be easily verified that
the  affine surface $K_1$ is of
logarithmic general type.

Assume that  the
Lang - Vojta conjecture
holds for $K_1$. It now follows from  lemma 1, that there is a uniform bound on
the number of integral points
on $E_t$: the integral points on $K_1$ are not Zariski dense, therefore the
integral points on $K_t$ are not Zariski dense, since they map to integral
points on $K_1$. Lemma 1 says that there is an open set $U\subset {\Bbb A}_1$
such that there is a  bound on  the number of points on $E_t$ for an
$S$-integer
$t$ in $U$; for the remaining finitely many integers $t$ we can use Siegel's
theorem. The same result can be obtained using any of the higher Kummer
varieties $(E\times\cdots\times  E)/(\pm 1)$.

Note that the integral points on $E_t$ are not the same as the
integral points on the N\'eron model. Suppose that $t$ is square free. Then a
N\'eron integral point $P$ is
integral on $E_t$, away from characteristic 2 and 3, if at a prime of additive
reduction $P$ does not reduce to the
component of the origin on the N\'eron model. In other words, even after
semistable reduction (obtained by taking $y' = \sqrt{t} y$), $P$ remains
integral on the N\'eron model. We call such
points {\em stably integral}.

{\bf Open problem:} Show that the $S$-integral points on $K_1$
are not Zariski dense. As a first step, describe the images of nontrivial
morphisms ${\Bbb
A}_1\setminus 0 \rar K_1$.

\section{Stably integral points}\label{stably}

\begin{dfn}
Let $E$ be an elliptic curve over a number field $K$, let $S$ be a finite set
of
places in $K$,
and $P$ be a $K$-rational point on $E$. We say that $P$ is {\em stably
$S$-integral}, written $P\in E(K,S)$  if the following holds: let $L$ be a
finite extension of $K$,
and let $T$ be the set of places above $S$, and assume that $E$ has semistable
reduction $\cle$ over $\co_{L,T}$; then  $P\in (\cle\setminus 0)(\co_{L,T})$.
In
other words, $P$ is
 integral on the  semistable model of $E\setminus 0$ over some
finite
field extension $L$ of $K$, where $T$ is the set of all places over $S$.
\end{dfn}

Stably integral points should be thought of as the rational  points which are
integral over the algebraic closure of the field. In this sense, they are a
good
analogue on elliptic curves, for rational points on curves of higher genus.

It is important to note that stably integral points can be described as the
integral points on a certain type of model of the curve.

\begin{dfn} Let $E$ be an elliptic curve over a number field $K$, let $S$ be a
set of places containing all places dividing 2 and 3, and let $\cle$
be the N\'eron model over $\oks$. Let $D_0$ be the zero section of $\cle$. Let
$S_a$ be the set
of places of additive reduction, and for a place $v$ let $\cle^0_v$ be the zero
component. Let $D= D_0\cup \bigcup_{v\in S_a} \cle^0_v$ and let $\cle_0 =
\cle\setminus
D$. We call $\cle_0$ the {\em Stably minimal model} of $E$.
\end{dfn}

\begin{prp} Let $S$ be a set of places containing all places dividing 2 and
3. Then the $S$-integral points on the stably minimal model are
precisely the stably $S$-integral points.
\end{prp}
{\bf Proof:}  One can prove this proposition using the explicit list of
possible reduction of the
N\'eron model and their semistable reduction (Tate's algorithm). If one
goes through this list, one
sees that the kernel of the semistable reduction map away from characteristic 2
and 3 is
precisely  the additive components of the identity on the N\'eron model. A
much more appealing proof follows directly from \cite{edix}, section 5
(especially remark 5.4.1): assume given
a field extension $L\supset K$ which is tamely and totally ramified at a given
prime $p$ (this can be assumed for a local field of semistable reduction of an
elliptic curves once
one avoids the primes dividing 2 and 3). Let $\cle_K, \cle_L$ be the N\'eron
models of $E$ over $K$ and $L$ respectively. In \cite{edix} one obtains a
description of
the map induced on N\'eron models $\cle_K\times\spec{\cal O}_L \rar \cle_L$,
and one there sees that the group of components of the reduction $(\cle_K)_p$
 maps isomorphically to the group of components of the fixed locus under
the Galois action of $(\cle_L)_p$. Therefore the kernel of the map of N\'eron
models is connected.

{\bf Remark:} in order to include primes over 2 and 3 one simply needs to
remove all
the additive components which are in the kernel of the semistable reduction
map. Since we are allowed to remove a finite number of places anyway, we can
ignore this problem altogether.

\section{The main theorem}\label{main}

{\bf Proof of the main theorem:} We may assume that $S$ contains all places
above 2 and 3, and that
$K$ contains $\zeta_3$.
 Let $X \rar U$ be the family of all elliptic curves given by the equation
$$y^2 = x^3 + ax + b,$$
over an open set $U$ in $\Bbb{A}^2$, with parameters $a,b$. Let $B_0$ be
any
irreducible closed subset in $U$ and $B$ a  compactification of $B_0$. We
can add level 3 structure to produce a semistable family $X_1\rar B_1$, over a
Galois, generically finite  cover $B_1$ of
$B$.  The  Galois group of the cover is some $G_1\subset Sl_2(\bff_3) =G$. We
have a natural map $X_1\rar E_1$,
coming from the moduli interpretation of $E_1$, which is $G_1$ equivariant.
We
have that  $X^n_{B_0}$ maps  to $((X_1)^n_{B_1})/G_1,$ which maps down to
$E_n/G$.
Since  the family $E_1$ is semistable, an $n$-tuple of stably integral
points on an elliptic curve gives rise to an integral point on
 $(E_n\setminus\Theta_n)/G$.

 Note that if $X_{B_0}\rar B_0$ is not isotrivial, then $X_{B_0}^n$ dominates
$E_n/G$.
Otherwise, its image is isomorphic to $(E\setminus 0)^n/Aut\,E$  for some fixed
elliptic curve
$E$, where $Aut\, E$ acts diagonally.

The following lemmas show that
$(E_n\setminus\Theta_n)/G$ is of logarithmic general type for large $n$, and
that for any
fixed elliptic curve $E$, $(E\setminus 0)^n/Aut\,E$ is of logarithmic general
type.
Assuming the
Lang-Vojta
conjecture, the integral points on the variety $((X_1)^n_{B_1}\setminus
\Theta)/G_1$ are not Zariski dense. By lemma 1, we obtain a uniform bound on
the
number of stably integral
points on all elliptic curves away from a closed subset $B'$ of $B$. By
Noetherian
induction, we have a bound on all elliptic curves. This gives the theorem.

\begin{lem}\label{flem}(Compare \cite{chm}, lemma 4.1) Let $X_0\subset X$ be an
open inclusion of an
irreducible
variety $X_0$ in a smooth complex projective irreducible $X$ of dimension $n$,
such that the complement
$D = X \setminus X_0$ is a divisor of normal crossings. Let $G$ be a finite
group acting on $ X,X_0,D$ compatibly. Let $\omega$ be a $G$- equivariant
logarithmic $k-$canonical form on $X_0$. If at any point $x$ of $X_0$ which is
fixed
by some element in $G$, the form $\omega$ vanishes to order at least
$C=k(|G|-1)$, then
$\omega$ descends to a regular logarithmic $k$-canonical form on any
desingularization of  $X_0/G$
\end{lem}

{\bf Proof:} let $Y_0$ be a desingularization of $X_0/G$ and let $Y$ be
a regular compactification, mapping to  $X/G$. Let $Z'$ be the graph of the
rational map $X\rar Y$. Let $Z$ be a $G$-equivariant desingularization of $Z'$,
and $Z_0$ the inverse image of $X_0$. Let $S$ be the branch locus of $Z$ over
$Y$. By a theorem of Hironaka, such
desingularizations may be chosen such that $(Z\setminus Z_0)\cup S$ is a
divisor
of normal crossings. Let $F_1$ be the closed set in $Y$ where the fibers in $Z$
are positive dimensional.   Let
$D_Z$ be the inverse image of $D$ in $Z$, $D_Y$ its image in $Y$. Let
$F_2\subset Y$  be the singular locus of $D_Y\cup S$. Clearly $F=F_1\cup F_2$
is of codimension at least 2 in $Y$. Note that away from $F_1$  the branch
locus $S$ is of codimension  1, since $Y$ is smooth.

Clearly $\omega$ descends to a  logarithmic form on $Y\setminus S$. It is
enough to show that it extends over $Y\setminus F$, since $F$ has codimension
at least 2.

Given a point   $y\in S\setminus F$ let $z\in Z$ be a point mapping to
it. We can choose formal coordinates $(z_1,z_2,\ldots,z_n)$ on $Z$ such that
$(y_1,z_2,\ldots,z_n)$ are coordinates on $Y$, with $y_1 = z_1^m$. Since we
removed the intersections of components of $D\cup S$, there are
only two cases to consider:

Case 1: $z\not\in D_Z$. We can write
$\omega =  f(z_1\ldots,z_m)z_1^C  (dz_1\wedge\cdots\wedge dz_m)^k$. We have
$dy_1 = m z_1^{m-1} dz_1$. Since $m<|G|$, we have that
$\omega = f z_1^{C-k(m-1)} (dy_1\wedge dz_2\wedge\cdots\wedge dz_m)^k$, is
regular, and since it is invariant it descends.

Case 2: $z\in D_Z$ and $z_1=0$ is the equation of $D_Z$. We can write
$\omega =  f(z_1\ldots,z_m)  (dz_1\wedge\cdots\wedge dz_m)^k/z_1^k$. Since
$m dz_1/z_1 =  dy_1/y_1,$ the invariance of $\omega$ means that  $f$
descends to $Y$, and
therefore
$\omega$ descends.

\begin{lem}\label{flem1}(Compare \cite{chm}, theorem 1.3)
\begin{enumerate} \item There exists a positive integer $n$ such that
$(E_n\setminus \Theta)/G$ ($G$ acting diagonally) is of
logarithmic general type.
\item For a fixed elliptic curve $E$, there is $n$ so that $(E\setminus
0)_n/Aut\,E$
($Aut\,E$ acting diagonally) is of logarithmic
general type.
\end{enumerate}
\end{lem}
{\bf Proof:} Let $S\subset E_1$ be any divisor containing the  locus of fixed
points of elements of $G$, and let $F$ be a fiber. Then the fixed
points in $E_n$ are contained in $S^n_{\bfp^1}$, the fibered product of $S$
with itself $n$ times. Recall that we have shown that
$L=\omega_{E/\bfp^1}(\Theta) $ is big; therefore for some large $k$, the
$\bfq$-line bundle $L(-(S+2F)/k)$ is big. This means that for $n=k|G|$, on
$E_n$ there are
many sections of $\omega^m_{E_n}$ vanishing to order $m|G|$ on the  fixed
points in $E_n$. As in \cite{chm}, lemma 2.1, it follows that there are also
many {\em invariant} sections
vanishing to such order. The proof of part (2) is identical.

\section{The uniform boundedness conjecture}\label{ubc}
It is well known that torsion points of high order on an elliptic curve are
integral; we will use this to study torsion points in terms
of integral points. As quoted in the introduction, a long standing conjecture
which was recently
proved by Merel  says that the order of a torsion point on an elliptic curve
over a number field is bounded in terms of the degree of the field of
definition only. We now indicate how Merel's theorem  follows from the
Lang - Vojta conjecture.

We start with the basic proposition which makes things work (see the case of
elliptic curves in \cite{oester}):

\begin{prp}
 Let $P$ be a torsion point on an abelian variety $A$,
both defined over a field $K$
of degree $d$. Denote by $n$ the order of $P$. Let $g$ be the dimension of
$A$ and let $C$ be the order of the group $Sp_g(Z/5Z)$. Assume that either $n$
is not a prime power, or $p^k = n$, such that $p^k-p^{k-1} > Cd$. Then
$P$ is stably integral on $A$, that is, its reduction at any place  on the
N\'eron model after semistable reduction is not the origin.
\end{prp}

{\bf Proof:} by adding level 3 structure we have (by a theorem of Raynaud) that
there is a field of degree at most $Cd$ where $A$ has semistable reductions
over all $p\neq 3$. Theorem IV.6.1 in \cite{silv} says that if $P$ is not
integral, then it is not integral at a place $\frak{p}$ above some prime $p$
where $n=p^k$;
and
the valuation satisfies $v(p)>p^k-p^{k-1}$. Here by definition,
$p=u\pi^{v(p)}$, where $u$ is a unit and $\pi$ a uniformizer of the valuation
ring. But $v(p)$ is
at most the degree of the field. We can similarly deal with primes over 3 by
adding level
5 structure instead. \qed

This following corollary is probably well known: torsion on abelian varieties
is bounded in terms of the degree of the field, the dimension and a prime of
potentially good reduction.

\begin{cor}
 For any triple $(d,g,p)$ there is an (explicit) integer $N$ such that if $K$
is a
number field of degree $d$, $A$ an abelian variety of dimension $g$ over $K$,
and $p$ is a rational prime over  which there is a place $\frak{p}$ of $K$
where $A$ has potentially good reduction, then $A(K)_{tors}<N$.
\end{cor}

{\bf Proof:} Let $L$ be a field of degree $\leq Cd$ over which $A$ has good
reduction at some prime $\frak{p}$ over $p$. Let $A(L)_{tors}\rar A_{\frak{p}}$
be the reduction map. By Weil's theorem, the image has cardinality $\leq
(1+p^{Cd/2})^{2g}$. But by the proposition, any point in the kernel is of
order $p^k$ satisfying $p^k-p^{k-1} < Cd$, which can be bounded as well.\qed

We would like to apply the correlation method to torsion points of high order
on elliptic curves, defined over all fields of degree $d$. By the proposition,
we may
use the fact that when $p$ is large these points are stably integral.

Since we want to show that there is a bound for torsion over number fields
depending only on the degree, we might as well assume that $E$ has level 3
structure: this has the effect of increasing the degree $d$ by a factor of 24.

 In \cite{kama}, Kamienny
and Mazur show that it is enough to
bound the order of prime torsion points. We will show the existence of a  bound
on prime order torsion points, assuming the Lang - Vojta conjecture.

Let $P$ be a torsion point of large prime order $p$ on an elliptic curve $E$
which has
level 3 structure, defined over some number field of degree $d$. The point $P$
gives a point on the surface $E_1$ introduced in \S\ref{level3}, defined over
the
same number field, and is in fact integral on $E_1\setminus \Theta$. The Galois
orbit of $P$ gives a $\bfq$-rational point on
the $d-$th symmetric power of $E_1$. We can do a bit better: fix an integer
$n$. Given $n$ torsion
points on an elliptic curve defined over the same number field (e.g. multiples
of a given torsion point) we in fact get a rational point on
$Y_n=\sym^d(E_n)$. In $Y_n$ there is a
divisor $\Theta_{Y_n}$ which consists of those tuples of  points such that
at least one point is the origin, and the points thus obtained are in fact
integral on the scheme $Y_n\setminus \Theta_{Y_n}$.

Given an auxiliary integer $k$, let $F_{n,k}$ be the Zariski closure in $Y_n$
of the set of all points
corresponding to Galois orbits of $n$-tuples of distinct torsion points of
prime order {\em larger than $k$} defined over fields of degree exactly $d$.

 By definition, if
$l>k$ then $F_{n,l}$ is contained in $F_{n,k}$. Let $F_n$ be the intersection
of
$F_{n,k}$ over all integers $k$. By the noetherian property of algebraic
varieties, $F_n = F_{n,k}$ for some $k$. What we want to show is that $F_n$
is empty. We will assume the contrary and derive a contradiction.

We have the natural symmetrization map $(E_n)^d \rar Y_n$. Let $G_n$ be the
inverse image of $F_n$ under this map.

We denote by $\pi^d_i:(E_n)^d \rar (E_1)^d$  the map induced from
$\pi_i:E_n\rar
E_1$.

The varieties $(E_n)^d$ can be viewed as compactified
semiabelian schemes over the  space $\bfp=(\bfp^1)^d$.

\begin{lem}
Let $G$ be a component of $G_n$.
There exists a closed subscheme $B\subset \bfp$, and  subvarieties $A_i \subset
E_1^d, 0< i<n+1$ mapping onto  $B$, such that the general fiber of $A_i$ over
$B$ is a
finite union of abelian subvarieties, and $G$ is a component of
the fibered product of $A_i$ over $B$. The varieties $A_i$ are not contained
in any diagonal in $E_1^d$ or in the theta divisor, nor in the locus of
singular
fibers.
\end{lem}

Proof: if $P$ is a torsion point of some prime order $p$ defined over a number
field, then
any multiple of it $kP$, for $k$ prime to $p$, is also torsion defined over the
same field. Fix an integer $1\leq i\leq n$. We look at the projection of
$q_i:G\rar E_{n-1}^d$ forgetting the $i$-th factor. It follows that each
fiber of
$G_n$ over $E_{n-1}^d$ is stable
under multiplication by $k$ for any integer $k$. We now use the trick of Neeman
and Hindry (see \cite{neeman} or \cite{hindry}), which tells us that a
subvariety of an abelian
variety which is stable under multiplication by all integers, is a union of
abelian subvarieties.

Let $G'\subset G_{n-1}$ be the image  of $G$ under $q_i$. For each point  $P\in
E_{n-1}^d$ in  $G'$ we have that
$q_i^{-1}(P)\cap F_n$ is a union of finitely many
 abelian
subvarieties  of $E_1^d$. Since a subvariety of a constant
abelian scheme is constant (say, by looking at torsion points), these abelian
subvarieties depend only on the image of $P$ in $\bfp$. Therefore there is a
subvariety $A_i$ as in the lemma such that $G$ is a component of
$(\pi^d_i)^{-1}A_i\cap q_i^{-1}G'$. By induction we obtain the product
structure.

Since the variety $G$ was obtained from the closure of Galois orbits of points
over fields of degree exactly $d$, none of them is fixed by any permutation,
and none is in the theta divisor. Similarly, we see that they are not contained
in the singular fibers of $E_n^d$.
\qed

We will now show that for high enough $n$, any candidate for a component of $F$
is of logarithmic general type. First note that, by noetherian induction, one
may assume that the base $B$ of $A_i$ remains constant as $n$ grows. We will
now see that if one of the $A_i$ appears many times in the product, then the
image $F$ of $G$ is of logarithmic general type.

\begin{lem}  Let
$B\subset
\bfp$ be an irreducible closed subvariety.
Let $A_i\subset (E_1)^d,\quad  1\leq i\leq m$ and $A_{m+1}$ be subschemes
mapping to
$B$ satisfying the conclusions in the previous lemma.
There is an integer $k_0$ such that for any $k>k_0$ and any $l_i\geq 0$ the
following
holds:

 Let $G$ be a component of the scheme
$$(A_1)^{l_1}_{B}\times_B\cdots\times_B (A_m)^{l_m}_B\times_B
(A_{m+1})^k_B.$$ Let $F$ be
the image of
$G$ in $Y_n$, where $n=l_1+\cdots +l_n+k$. Let $F' = F\setminus\Theta_{Y_n}$.
Then $F'$ is of logarithmic general  type.
\end{lem}

{\bf Proof:}  Let $M_i$ be the dimension of the fibers of $A_i$ over $B$. For
each
choice of a subset $J_i$ of $\{1,\ldots,d\}$ of size $M_i$, we have a variety
$E_{J_i,B}$,
the
pullback of
$(E_1)^{M_i}$ to $B$ along the projection $\pi_{J_i}$ to the factors in $J_i$.
Since
$A_i$ is
not contained in the theta divisor of $E_1^d$, we have that $A_i$ surjects
generically finitely onto
$E_{J_i,B}$, whenever $J_i$ has size $M_i$. We treat $A_{m+1}$ a bit
differently:  using $d$ different generically finite  surjections $A_{m+1}\rar
E_{J'_i,B}$ where $J'_i=\{i,\ldots,(i+M_{m+1} \mod d)\}\subset \{1,\ldots,d\}$,
we can cook up a special generically finite surjection: write $k=qd+r$, then
we map map $(A_{m+1})^k_B\rar (E_{J_{m+1},B})^r_B  \times_B (E_{qm})^d|_B $.

In order to deal with the singularities, we desingularize the base: $B'\rar B$.
Now the pullback of the product of $E_{J_i,B}$ to $B'$  has  semistable
fibers,  therefore has log canonical singularities as in lemma 2. Choose a
canonical divisor
$K_{B'}$. Choose an effective divisor
$H\subset B$ such that the pull-back of $H$ to $B'$ is bigger than $-K_{B'}$.
If $G'$ is a
desingularization of $G$, it admits a generically finite surjection
$$p:G\rar V= (E_{J_1,B})^{l_1}_B \times_B\cdots\times_B
(E_{J_m,B})^{l_m}_B\times_B
(E_{J_{m+1},B})^r_B  \times (E_{qm})^d|_B.$$
Notice that $V$ is the restriction to $B$ of a variety of the form
$E_{r_1}\times \cdots \times E_{r_d}$, and if $k$ is large then {\em each of
the $r_i$} is large as well.
We can choose $G'$ so  that it maps to $B'$. Let $V'$ be the pullback of $V$ to
$B'$.

We wish to use the sections of powers of the logarithmic relative dualizing
sheaf of the product variety $V$  to construct
differential forms on $F$.  In view
of  lemmas \ref{flem},\ref{flem1}  we need the sections to vanish sufficiently
along the preimage of $H$,  and their pullback to $G$ should vanish along the
fixed points $\Delta$ of the symmetric group action to sufficiently high order.

Since the relative dualizing sheaf $\omega_{E_1/\bfp^1}(\Theta)$ is nef and
big, then the sheaf $\omega_{E_{J,B}/B}(\Theta)$ is nef, and the sheaf
$\omega_{(E_{1})^d|_B/B}(\Theta)$ is
nef and big. We have
an injection $p^*(\omega^m_{V'/B'})(-mH+m\Theta) \rar
\omega^m_{G'}(m\Theta)$, since  $V'\setminus \Theta$ has log
canonical singularities. Since each of the $r_i$ in the description of $V$ can
be made as large as we wish, the argument of lemma \ref{flem1} shows that $F$
is of logarithmic general type. \qed

We can now show by induction on $M$ that the relative dimension of $G_{n+1}$
over
$G_{n}$ is at least $M+1$, thus obtaining a contradiction. Clearly the
relative dimension is at least 1. If for some $n$ the relative dimension is
precisely $M$, then by induction, using the embedding $G_{n+k}\subset
(G_{n+1})^k_{G_n}$, the relative
dimension of $G_{n+k}$ over $G_n$ is $Mk$, and therefore there is a component
$G$ of  $G_{n+k}$ of relative dimension $Mk$. From lemma 6 it follows that $G$
is a component of a product variety of the form described in lemma 7, and
therefore $F\setminus \Theta$ is of logarithmic general type. The Lang - Vojta
conjecture implies that the integral points on $F$ are not dense, contradicting
the definition of $F$.

We arrived at a contradiction, therefore $F_n$ must be empty, and we  conclude
that there is a bound for torsion points of prime order. \qed

\hspace{8cm}\parbox{8cm}{{\sc Dan Abramovich}\\ Department of Mathematics\\
Boston University\\
111 Cummington Street
Boston, MA.  02215 \\
U.S.A.\\ e-mail:
abrmovic@@math.bu.edu\\
 \today}

\end{document}